\documentclass[12pt,a4paper]{amsart}
\usepackage[english]{babel}
\usepackage{a4wide,amssymb,xypic}
\newtheorem{thm}{Theorem}[section]
\newtheorem{corol}[thm]{Corollary}
\newtheorem{lemma}[thm]{Lemma}
\newtheorem{prop}[thm]{Proposition}
\newtheorem{defin}[thm]{Definition}

\theoremstyle{remark}
\newtheorem{rem}[thm]{Remark}
\newtheorem{ex}[thm]{Example}

\newenvironment{example}{\begin{ex}\rm}{\qee\end{ex}}
\newcommand{\qee}{\mbox{\hspace{0.2mm}}\hfill$\triangle$}

\newcommand{\Oc}{\mathcal O}

\newcommand{\PP}{{\mathbb P}}

\newcommand{\R}{\mathbb R}

\newcommand{\Q}{\mathbb Q}

\def\rk{\operatorname{rk}}
\def\coker{\operatorname{coker}}

\def\dim{\operatorname{dim}}

\newcommand\grass{\mbox{Gr}}
\newcommand\hgrass{{\mathfrak{Gr}}}
\newcommand{\cO}{{\mathcal O}}
\newcommand{\OPE}{{\cO_{\PP E}}}
\newcommand{\fE}{{\mathfrak E}}
\newcommand{\fF}{{\mathfrak F}}
\newcommand{\fG}{{\mathfrak G}}
\newcommand{\fQ}{{\mathfrak Q}}

\newcommand{\OPQ}[1]{{\mathcal O_{\PP Q_{#1}}}}
\begin{document}
\rightline{SISSA Preprint 39/2005/fm} \bigskip\bigskip
\title[Numerically flat Higgs vector bundles]{Numerically flat Higgs vector bundles}
\date{16 June 2005, revised 17 March 2006}
\thanks{This research was partly supported by the Spanish {\scriptsize
DGES} through the research project BFM2003-00097, by ``Junta de Castilla y Le{\'o}n''
through the research project SA114/04, by Istituto Nazionale per l'Alta Matematica,  and by the Italian National Research Project ``Geometria delle variet\`a algebriche''.}
 \maketitle \thispagestyle{empty}
\begin{center}{\sc Ugo Bruzzo} \\[2pt]
Scuola Internazionale Superiore di Studi Avanzati,\\ Via Beirut 2-4, 34013
Trieste, Italia;\\[2pt] Istituto Nazionale di Fisica Nucleare, Sezione di Trieste \\[2pt]  E-mail {\tt bruzzo@sissa.it} \\[10pt]
{\sc Beatriz Gra\~na Otero} \\[2pt]
Departamento de Matem\'aticas and Instituto de F\'\i sica \\ Fundamental y Matem\'aticas, Universidad de Salamanca, \\ Universidad de Salamanca,
\\  Plaza de la Merced 1-4, 37008 Salamanca, Espa\~na\\[2pt] E-mail {\tt beagra@usal.es}
\end{center}
\bigskip
\begin{abstract} After providing a suitable definition of  numerical effectiveness for Higgs bundles, and a related notion of numerical flatness,  in this paper we prove, together with some side results, that all Chern classes of a Higgs-numerically flat Higgs bundle vanish,
and that a Higgs bundle is Higgs-numerically flat if and only if it is has
a filtration whose quotients are  flat stable Higgs bundles. We also study the relation between
these numerical properties of Higgs bundles and (semi)stability.
\end{abstract}

\newpage

\section{Introduction.}
A notion of ampleness for a vector bundle $E$ on a scheme $X$ was
introduced in \cite{H66}: $E$ is said to be ample if the tautological
bundle $\OPE(1)$ on the projectivization $\PP E$ is ample.
A related notion is that of numerical effectiveness:
$E$ is numerical effective if $\OPE(1)$ is so (in some literature,
the terms ``pseudo-ample'' or ``semi-positive''  are used instead of ``numerically effective'').
Moreover, $E$ is said to be \emph{numerically flat} if both $E$
and the dual bundle $E^\ast$ are numerically effective.

Properties of numerically effective and numerically flat vector
bundles were studied in several papers, see e.g.
\cite{Fu83,CP91,DPS94}. In particular, all Chern classes of a
numerically flat vector bundle vanish (the same result has been
proved for principal bundles in  \cite{BS05}). In this paper we
explore some properties of  the Higgs bundles that satisfy
numerical effectiveness or flatness conditions; these are a
modification of those  introduced in \cite{BHR05}. We give a
recursive definition: a Higgs bundle $\fE$ of rank one is said to
be \emph{Higgs-numerically effective} (resp.~\emph{Higgs-ample})
if it is numerically effective (resp.~ample) in the usual sense,
while if its rank is greater than one we require that its
determinant bundle is nef (resp.~ample), and all universal Higgs
quotient bundles on the Grassmannian varieties of Higgs quotients
are Higgs-numerically effective (resp.~Higgs-ample). Finally, a
Higgs bundle $\fE$ is said to be \emph{Higgs-numerically flat} if
both $\fE$ and its dual are Higgs-numerically effective.

Among other things, we prove that the Chern classes of a
numerically flat Higgs bundle vanish. One encounters here a
similar situation as in the case of the Bogomolov inequality:
semistable (ordinary) bundles satisfy the Bogomolov inequality,
but this is also implied by the weaker property of
Higgs-semistability. In the same line, the condition of
Higgs-numerical flatness is weaker than ordinary numerical
flatness, but it is nevertheless enough to ensure that all Chern
classes vanish. We also prove that  a Higgs bundle is
Higgs-numerically flat  if and only if admits a filtration whose
quotients are   flat stable Higgs bundles.

One of the main ingredients of our proofs is the numerical characterization
of the semistability of Higgs-bundles that has been established in
\cite{BHR05} (see Theorem \ref{1.3BHR} in this paper).

{\bf Acknowledgment.} The authors thank G.~Ottaviani for a useful remark.

\section{Ample and numerically effective Higgs bundles}

All varieties are projective varieties
over the complex field.
Let $E$ be a vector bundle of rank $r$ on $X$, and let $s$ be a positive
integer less than $r$. We shall denote by $\grass_s(E)$  the Grassmann bundle of $s$-planes in
$E$, with projection $p_s : \grass_s(E) \to X$. When
$s=1$, the Grassmann bundle reduces to
the projectivization of $\PP E$ of  $E$, defined as
$\PP E = \mbox{\bf Proj}(\mathcal S(E))$, where $\mathcal
S(E)$ is the symmetric algebra of the sheaf of sections of $E$.

The Grassmann bundle $\grass_s(E)$ is a parametrization of the rank $s$ locally-free quotients of $E$. There
is a universal exact sequence
\begin{equation}\label{univ}
0 \to S_{r-s,E} \xrightarrow{\psi} p_s^*(E) \xrightarrow{\eta} Q_{s,E}
\to 0
\end{equation}
of vector bundles on $\grass_s(E)$, with $S_{r-s,E}$ the universal rank $r-s$ subbundle and $Q_{s,E}$ the universal rank $s$ quotient bundle.

\begin{defin} A Higgs sheaf $\fE$ on $X$ is a coherent sheaf $E$ on $X$ endowed with
a morphism $\phi \colon E \to E \otimes \Omega_X$ of $\Oc_X$-modules such that
$\phi\wedge\phi=0$, where $\Omega_X$ is the cotangent sheaf to $X$. A Higgs subsheaf
$F$ of a Higgs sheaf $\fE=(E,\phi)$ is a subsheaf of $E$ such that $\phi(F)\subset
F\otimes\Omega_X$. A  Higgs bundle is a Higgs sheaf $\fE $ such that $E$ is a
locally-free $\Oc_X$-module.
\end{defin}

There exists a stability condition for Higgs sheaves analogous to that for ordinary
sheaves, which makes reference only to $\phi$-invariant subsheaves.

\begin{defin} Let $X$ be a smooth projective variety equipped with a
polarization. A Higgs sheaf $\fE=(E,\phi) $ is semistable (resp.~stable) if $E$ is
torsion-free, and $\mu(F)\le \mu(E)$ (resp. $\mu(F)< \mu(E)$) for every proper
nontrivial Higgs subsheaf $F$ of $\fE$.\end{defin}

Given a Higgs bundle $\fE $, we may construct closed subschemes
$\hgrass_s(\fE)\subset \grass_s(E)$ pa\-ram\-e\-ter\-iz\-ing rank
$s$ locally-free Higgs quotients, i.e., locally-free quotients of
$E$ whose corresponding kernels are $\phi$-invariant. With
reference to the exact sequence eq.~\eqref{univ}, we define
$\hgrass_s(\fE)$ as the closed subvariety of $\grass_s(E)$ where
the composed morphism
$$
(\eta\otimes1)\circ p_s^\ast(\phi) \circ \psi\colon S_{r-s,E}\to Q_{s,E}\otimes
 p_s^\ast\Omega_X
$$
vanishes. We denote by $\rho_s$ the projections $\hgrass_s(\fE)\to X$. The
restriction of \eqref{univ} to the scheme $\hgrass_s(\fE)$ provides the universal
exact sequence
$$
0\to S_{r-s,\fE}\xrightarrow{\psi} \rho_s^\ast (E)\xrightarrow{\eta}
 Q_{s,\fE}\to 0
$$
and $Q_{s,\fE}$ is a rank $s$ universal Higgs  quotient vector bundle, i.e., for every
morphism $f\colon Y\to X$ and every  rank $s$ Higgs quotient $F$ of $f^\ast E$ there
is a morphism $\psi_F\colon Y\to \hgrass_s(\fE)$ such that $f=\rho_s\circ \psi_F$
and $F\simeq\psi_F^\ast (Q_{s,\fE})$. Note that the kernel $S_{r-s,\fE}$ of the morphism
$\rho_s^\ast (E )\to Q_{s,\fE}$ is $\phi$-invariant.

The scheme $\hgrass_s(\fE)$ will be called  the \emph{Grassmannian of locally free rank $s$ Higgs quotients} of   $\fE$.

We give now our definition of ampleness and numerical
effectiveness for Higgs bundles. This is inspired by, but is
different from, the definition given in \cite{BHR05} (Def.~4.1):
on the one hand it is less restrictive (though it is enough for
proving the results we have in mind), while on the other hand we
add a requirement on the determinant bundle to avoid the inclusion
of somehow pathological situations, cf.~Example \ref{patho} below.
\begin{defin}  \label{moddef} A Higgs bundle $\fE$ of rank one is said to be Higgs-ample
if it is ample in the usual sense. If $\rk \fE \geq 2$ we require
that:
\begin{enumerate} \item  all bundles $Q_{s,\fE}$ are Higgs-ample;
\item the line bundle $\det(E)$ is ample.
\end{enumerate}
A Higgs bundle $\fE$ of rank one is said to be Higgs-numerically
effective if it is numerically effective in the usual sense. If
$\rk \fE \geq 2$ we require that:
\begin{enumerate} \item all bundles $Q_{s,\fE}$ are Higgs-nef; \item the line bundle $\det(E)$ is nef.
\end{enumerate}
If both $\fE$ and $\fE^\ast$ are Higgs-numerically effective, $\fE$ is said to be Higgs-numerically flat.
\end{defin}

For short we shall use the abbreviations \emph{H-ample} for
Higgs-ample,  \emph{H-nef} for Higgs-numerically effective, and
\emph{H-nflat} for Higgs-numerically flat. Note that if
$\fE=(E,\phi)$, with $E$ ample in the usual sense (resp.~nef) in
the usual sense, than $\fE$ is H-ample (resp.~H-nef). Moreover, if
$\phi=0$, the Higgs bundle $\fE=(E,0)$ is H-ample (resp.~H-nef) if
and only if $E$ is ample (resp.~nef) in the usual sense.

\begin{example} Examples of Higgs bundles that are H-nflat but not numerically flat
as   ordinary bundles may be constructed in terms of a Higgs
bundle $\fE=(E,\phi)$ which is semistable as a Higgs bundle but
not as an ordinary bundle. Let $\fF=\fE\otimes\fE^\ast=(E\otimes
E^\ast,\psi)$. Since $\fE$ is semistable, and $c_1(E\otimes
E^\ast)=0$, by Proposition \ref{box2} below the Higgs bundle $\fE$
is H-nef, and hence H-nflat. On the other hand if  $E\otimes
E^\ast$ were numerically flat it would be semistable as a vector
bundle (cf.~Corollary \ref{corss}), and then $E$ would be
semistable as well --- a case we are excluding.
\end{example}

\begin{example}\label{patho} Let us motivate the appearance of the condition that $\det(E)$ is ample, or nef, in Definition \ref{moddef}. Let $X$ be a smooth projective curve, and let $\fE=(E,\phi)$ be a rank 2 nilpotent Higgs bundle, i.e.,
$E=L_1\oplus L_2$ (where $L_1$, $L_2$ are line bundles),
and $\phi\colon L_1\to L_2\otimes \Omega_X$, $\phi(L_2)=0$. It is shown
in \cite{BHR05} that $\hgrass_1(\fE) = \PP(L_1) \cup \PP (Q)$, where
$$ Q = \coker (\phi\otimes\mbox{id}) \colon E\otimes T_X \to E\,,$$
and $T_X$ is the tangent bundle to $X$. This implies that $\fE$
has only  two Higgs quotients, i.e., $L_1$ and $\bar Q$ which is $Q$ modulo torsion.
Note that $\deg(\bar Q) \ge \deg(L_1)$. If the genus of $X$ is at least 2,
one can for instance take $\deg(L_1)=0$ and $\deg(L_2) =-2$. Without the
  condition that $\det(E)$ is nef this Higgs bundle, which has negative degree, would be H-nef.  \end{example}

We prove now some properties of H-nef Higgs bundles that will be useful in the sequel. These generalize properties given in \cite{H66,CP91} for ordinary vector bundles.

\begin{prop} \label{properties}
Let $X$ be a smooth projective variety.
\begin{enumerate}
    \item A Higgs bundle $\fE=(E,\phi)$ on $X$ is H-nef if and only if
    the Higgs bundle $\fE\otimes \cO_X(D) = (E\otimes \cO_X(D), \phi\otimes \mbox{id})$
    is H-ample for every ample $\Q$-divisor $D$ in $X$.
    \item If $f\colon Y \to X$ is a finite surjective morphism of smooth projective
    varieties, and $\fE$ is a Higgs bundle on $X$, then $\fE$ is H-ample (resp.~H-nef) if and only if $f^\ast\fE$ is H-ample (resp.~H-nef).
     \item Every quotient Higgs bundle of a H-nef Higgs bundle
    $\fE$ on $X$ is H-nef.
\end{enumerate}
\end{prop}

\begin{proof} (i) If $\rk \fE = 1$ this is a standard property of nef line bundles.
Then we use induction on $\rk \fE$ using the fact that under the
isomorphism $\hgrass_s(\fE\otimes \cO_X(D)) \simeq \hgrass_s(\fE)$
the universal quotients of $\fE\otimes \cO_X(D)$ are identified
with $Q_{s,\fE} \otimes \rho_s^\ast\cO_X(D)$.

(ii) Again, this is standard in the rank one case \cite{H66}. In
the higher rank case we first notice that $\det(f^\ast(E))\simeq
f^\ast(\det(E))$, so that the condition on the determinant is
fulfilled. Moreover, by functoriality the morphism $f$ induces a
morphism $\bar f \colon \hgrass_s(f^\ast\fE)\to \hgrass_s(\fE)$,
and $Q_{s,f^\ast\fE} \simeq \bar f^\ast (Q_{s,\fE})$. One
concludes by induction.

(iii) Let $\fF=(F,\phi_F)$ be a rank $s$ Higgs quotient of $\fE$.
This corresponds to a section $\sigma\colon X \to \hgrass_s(\fE)$
such that $\fF \simeq \sigma^\ast(Q_{s,\fE})$. Since $Q_{s,\fE}$
is H-nef, $\fF$ is H-nef as well by the previous point.
\end{proof}

In \cite{Mi} Miyaoka introduced a numerical class $\lambda$ in the projectivization
$\PP E  $ which, when $X$ is a curve, is nef if and only if $E$ is semistable. In
\cite{BHR05} some  generalizations of the class $\lambda$ were introduced. In the
case of ordinary bundles one defines, for $s=1,\dots,r-1$,
$$\lambda_{s,E}= \left [ c_1(\OPQ{s,E}(1)) \right ] -\frac1r \; q_s^\ast ( c_1(E))
\in N^1(\PP Q_{s,E})\,,$$
where $q_s\colon \PP Q_{s,E} \to X$ is    the
natural epimorphism, and
$$\theta_{s, E}=[c_1(Q_{s,E})]-\frac sr\, p_s^\ast (c_1(E)) \in
N^1(\hgrass_s(E))\,.$$
Here, for every scheme $Z$, we denote by $N^1(Z)$ the
vector space of divisors modulo numerical equivalence: $$N^1(Z) =
\frac{\mbox{Pic}(Z)}{num.eq.} \otimes \R.\,.$$

It will be useful to have a formula relating the first Chern class
of the universal bundles to the classes $\theta_{s,E}$, that
we write in the following form.

\begin{lemma} \label{phidias} One has
$$c_1(S_{r-s,E}) = -\theta_{s,E} +\frac{r-s}{r}\,p_s^\ast(c_1(E))\,.$$
\end{lemma}
\begin{proof} The result is obtained by tensoring
the dual of the exact sequence \eqref{univ} by $Q_{s,E}$ and
recalling that $T_{\operatorname{Gr}_s(E)/X}\simeq
S_{r-s,E}^\ast\otimes Q_{s,E}$, cf.~\cite{Fu84}.
\end{proof}

In the case of a  Higgs bundle $\fE$ on a smooth projective variety $X$ one
defines
$$\lambda_{s,\fE}= \left [ c_1(\OPQ{s,\fE}(1)) \right ] -\frac1r \pi^\ast_s ( c_1(E))
\in N^1(\PP Q_{s,\fE})$$
$$\theta_{s, \fE}=[c_1(Q_{s,\fE})]-\frac sr\, \rho_s^\ast (c_1(E)) \in
N^1(\hgrass_s(\fE)),$$ where
$\pi_s\colon \PP Q_{s,\fE}  \to X$   and
 $\rho_s \colon \hgrass_s(\fE)\to X$ are again the
natural epimorphisms.

Miyaoka's criterion for semistability has been generalized in \cite{BHR05} to Higgs
bundles on smooth projective varieties of any dimension (the same criterion has been generalized to principal bundles in \cite{BB04}). Let $\Delta(E)$ be the
characteristic class
$$\Delta(E) =  c_2(E)-\frac{r-1}{2r}c_1(E)^2 = \frac{1}{2r}c_2(E\otimes E^\ast)\, .$$
Following theorem was partly proved in \cite{BHR05}.
\begin{thm} \label{1.3BHR}
Let $\fE$ be a Higgs bundle on a smooth projective variety. The following conditions
are equivalent.
\begin{enumerate}
    \item All classes $\lambda_{s, \fE}$ are nef, for $0<s<r$.
    \item $\fE$ is semistable and $\Delta(E)=0$.
    \item All classes $\theta_{s, \fE}$ are nef, for $0<s<r$.
    \item For any smooth projective curve $C$ in $X$, the restriction $\fE_{\vert C}$ is semistable.
\end{enumerate}
\end{thm}
\begin{proof}
The equivalence between (i) and (ii) was proved in \cite{BHR05}.

(i) implies (iii): the class $\lambda_{s,\fE}$ may be regarded as
the numerical class of the hyperplane bundle of the Higgs
$\Q$-bundle $\fF_s = \fQ_{s, \fE} \otimes
\rho_s^{\ast}(\det^{-1/r}(\fE))$, which therefore is nef. As a
consequence, the class $\theta_{s,\fE}=[c_1(\fF_s)]$ is nef.

(iii) implies (iv): the restriction $\theta_{s,\fE_{\vert C}}$ of
$\theta_{s, \fE}$ to $\hgrass_s(\fE_{\vert C})$ is nef. If $\fE'$
is a rank $s$ Higgs quotient of $\fE_{\vert C}$, let $\sigma\colon
C \to\hgrass_s(\fE_{\vert C})$ be the corresponding section. Then
$$\theta_{s,\fE_{\vert C}}\cdot [\sigma(C)] =
s(\mu(E')-\mu(E_{\vert C}))\ge 0 $$ so that $\fE_{\vert C}$ is
semistable.

(iv) implies (i): Let $C'$ be an irreducible curve in $\PP
Q_{s,\fE}$ and let $\overline{f} \colon C \to C'$ be it is a
normalization. If $f = \pi_s \circ \bar{f}$, the pullback
$f^{\ast} \fE$ is semistable, so that by Miyaoka's criterion the
divisor $\lambda_{s, f^{\ast} \fE}$ is nef. Since $$\lambda_{s,
\fE} \cdot [C'] = \deg(\lambda_{s, f^{\ast} \fE}) \geq 0$$ the
divisors $\lambda_{s,\fE}$ are nef.
\end{proof}
This admits a simple corollary \cite{BHR05}.
\begin{corol} \label{corral}
A semistable Higgs bundle $\fE=(E,\phi)$ on an $n$-dimensional polarized
smooth projective variety $(X,H)$ such that $c_1(E)\cdot H^{n-1}=\mbox{\rm
ch}_2(E)\cdot H^{n-2} = 0$ is H-nflat.
\end{corol}

Theorem \ref{1.3BHR} makes use of Theorem 2 in \cite{S92}, which
will also be further needed in the present paper. We recall it
here in a simplified form which is enough for our purposes.

\begin{thm} Let $\fE=(E,\phi)$ be a semistable Higgs bundle on an $n$-dimensional polarized smooth projective variety $(X,H)$, and assume
$c_1(E)\cdot H^{n-1} = \operatorname{ch}_2(E)\cdot H^{n-2} = 0$.
Then $\fE$ admits a filtration whose quotients   are stable and
have vanishing Chern classes. \label{thm2} \end{thm}

\section{Numerically flat Higgs bundles and stability}

In this section we explore some properties of Higgs bundles
related to the notion of H-numerical effectiveness and H-numerical flatness. We start by showing that all Chern classes of an H-nflat Higgs bundle
vanish.

\begin{thm} Let $\fE=(E, \phi)$ be an H-nef Higgs bundle on a smooth polarized projective variety $(X,H)$ whose first Chern class is numerically equivalent to zero,
$c_1(E) \equiv 0$. Then all Chern classes  $c_k(E)$ vanish. \label{firstday}
\end{thm}
\begin{proof} Since all bundles $Q_{s,\fE}$ are H-nef, and
$c_1(E) \equiv 0$, the classes $\theta_{s,\fE}$ are nef, so that
$\fE$ is semistable and $\Delta(E)=0$ by Theorem \ref{1.3BHR}.
Since $c_1(E)\equiv 0$ then $\mbox{ch}_2(E)\cdot H^{n-2}=0$ (where
$n=\dim X$), and Theorem \ref{thm2} implies that all Chern classes
$c_k(E)$ vanish.    \end{proof}

\begin{corol} Let $\fE=(E, \phi)$ be a Higgs bundle on a smooth polarized projective variety $X$.
If $\fE$ is   H-nflat, then all Chern classes  $c_k(E)$ vanish.
\end{corol}

\begin{proof}  Since $\det(E)$   is numerically flat,    the
class $c_1(E)$ is numerically equivalent to zero. Moreover, $\fE$ is H-nef,
so that Theorem \ref{firstday} applies. \end{proof}

The next result generalizes Corollary 3.6 in \cite{BHR05} and Theorem 1.2 in
\cite{G79}. The proof does not differ much from the one given in \cite{G79} but we
include it here for the reader's convenience.

\begin{prop} Let $\fE=(E,\phi)$ be a  Higgs bundle on a smooth
projective variety $X$ such that all classes $\lambda_{s,\fE}$ are nef.
\begin{enumerate} \item  If the  class $c_1(E)$ is nef, then all universal
quotient bundles $Q_{s,\fE}$ are nef (so that $\fE$ is H-nef).
\item  If $X$ is a curve and $c_1(E)$ is ample, then all universal
quotient bundles $Q_{s,\fE}$ are ample (so that $\fE$ is H-ample).
\item If  $c_1(E)$ is positive (i.e., $c_1(E)\cdot[C]>0$ for all
irreducible curves $C\subset X$), then the class $c_1(Q_{s,\fE})$
is positive for all $s$.
\end{enumerate} \label{box2}
\end{prop}
\begin{proof} (i).  If $Q_{s,\fE}$ is not nef there is an irreducible curve $C\subset \PP Q_{s,\fE}$ such that $c_1(\cO_{\PP
Q_{s,\fE}}(1))\cdot [C]<0$. Let $f\colon C'\to C $ be the normalization of $C$, and let
$p\colon C'\to\hgrass_s(\fE)$ be the induced map. If $L$ is the pullback of $\cO_{\PP
Q_{s,\fE}}(1)$ to $C'$, then $L$ is a Higgs quotient of $p^\ast\circ \rho_s^\ast (E)$,
and
$$ \deg(L) = [f(C')]\cdot c_1(\cO_{\PP Q_{s,\fE}}(1))<0\,.$$
On the other hand, one has
$$\deg(p^\ast\circ \rho_s^\ast (E)) = [p(C')]\cdot c_1(\rho_s^\ast (E))\ge 0$$
since $c_1(E)$ is nef, so that
\begin{equation}\label{contr}
\mu(L) < \mu(p^\ast\circ\rho_s^\ast (E)).
\end{equation}
Now, in view of Theorem   \ref{1.3BHR}, the fact that all classes
$\lambda_{s,\fE}$ are nef implies that $\fE$ is semistable, and also that the
restriction of $\fE$ to any smooth projective curve in $X$ is semistable. Combining
this with Lemma 3.3 in \cite{BHR05}, one shows that $p^\ast\circ\rho_s^\ast (E)$
is semistable. But then eq.~\eqref{contr} is a contradiction.

(ii). This proof is a slight variation of the previous one, due to
the fact that Nakai's criterion for ampleness requires to check
positive intersections with subvarieties of all dimensions. Let
$C$ be a smooth projective curve and $f\colon C\to X$ a morphism
which is of degree larger than $r=\rk E$. Given a point $p\in C$
let $F$ be the class of the fibre of $\PP(f^\ast Q_{s,\fE})$ over
$p$. The Higgs bundle $\fE'=f^\ast\fE\otimes \cO_C(-p)$ is
semistable by the same argument as in the previous proof.
Moreover, $\deg (\fE')>0$, so that $\fE'$ is H-nef by the previous
point. If $L$ is the pullback to $C$ of the bundle $\cO_{\PP
Q_{s,\fE}}(1)$, then $L(-F)$ is nef since it is the hyperplane
bundle in $\PP Q_{s,\fE'}$. If $V$ is any subvariety of
$\PP(f^\ast Q_{s,\fE})$ of dimension $k$, then
$c_1(L)^k\cdot[V]>0$, so that $L$ is ample. Thus the pullback of
$\fE$ to $C$ is H-ample, and hence $\fE$ is H-ample as well by
Proposition \ref{properties}.

Claim (iii) is proved as Claim (ii).
\end{proof}

\begin{corol} Given a Higgs bundle $\fE=(E,\phi)$, if
all classes $\lambda_{s,\fE}$ are nef, and $c_1(E)$ is numerically equivalent
to zero, then $\fE$ is H-nflat.
\end{corol}

We study now the relations between the conditions of H-numerical
effectiveness and flatness and (semi)stability.

\begin{prop}\label{ss}
Let $\fE=(E,\phi)$ be a  Higgs bundle on a smooth polarized projective variety $X$,
such that all universal quotients $Q_{s,\fE}$ and $Q_{s,\fE^\ast}$ are nef. Then $\fE$ is semistable. If $\deg(E)\ne 0$, then $\fE$ is stable.
\end{prop}
\begin{proof} Under the isomorphism $\hgrass_{r-s}(\fE^\ast)\simeq \hgrass_s(\fE)$
the bundle $S_{r-s,\fE}$  is identified with $Q_{r-s,\fE^\ast}^\ast$. Therefore all the universal quotient bundles $Q_{s,\fE}$ and the
bundles $S_{r-s,\fE}^\ast$ on $\hgrass_s(\fE)$ are nef.
From Lemma \ref{phidias} we have, after restricting to
$\hgrass_s{\fE}$,
\begin{equation}\label{c1nef} c_1(S_{r-s,\fE}^\ast) = \theta_{s,\fE} +  \frac{s-r}{r} {\rho_s}^\ast (c_1(E))\,.\end{equation}
By   \cite[Prop. 1.2 (11)]{CP91}  this class is nef.

Assume for a while that $X$ is a curve, and let us suppose that
$\deg(E) \ge 0$. By a slight generalization of \cite[Prop.
1.8(i)]{DPS94} or \cite[Prop. 2.2]{Fu83}, the class $p_s^\ast
(c_1(E))$ is positive, and as $\hgrass_s(\fE)$ is a closed
subscheme of $\grass_s(E)$, the class $\rho_s^\ast (c_1(E))$ is
positive as well.  But since $c_1(S_{r-s,\fE}^\ast)$ is nef this
implies that all classes $\theta_{s, \fE}$ are nef and so from
Proposition \ref{1.3BHR} it follows that $\fE$ is   semistable.

If $\deg(E) \le 0$, the same argument shows that $\fE^\ast$ is semistable,
and then $\fE$ is semistable as well.

We now show that if $\deg(E)\ne 0$ then $\fE$ is stable. Assume for instance
that $\deg(E)>0$. Proposition \ref{box2} proves that in this case $c_1(Q_{s, \fE})>0$ for all $s$. Without loss of generality we may assume
that $\hgrass_s(\fE)$ has a   section $\sigma : X \to
\hgrass_s(\fE)$. Then the bundle $Q_\sigma = \sigma^\ast(Q_{s, \fE})$ is  an ample Higgs
quotient of $\fE$. So one has   the   exact sequence
$$ 0 \to K \to E \to Q_\sigma \to0$$
 and $ -c_1(K) = c_1(Q_\sigma \otimes \det^{-1}E) = \sigma^\ast(c_1(S_{r-s,\fE}^\ast))$
 is nef as well. Thus $c_1(K) \leq 0$ and $\mu(K) < \mu(E)$. Hence $\fE$ is  stable.
If $\deg(E)<0$ by applying the same argument to the dual of $\fE$ we obtain that
$\fE^\ast$ is stable, and hence $\fE$ is  stable again.

These results  are then extended to higher dimensional $X$ with the usual induction on the dimension of $X$, by considering a smooth divisor in the linear system $\vert mH\vert$ for $m$ big enough.
\end{proof}

\begin{corol} An H-nflat Higgs bundle is semistable. \label{corss}
\end{corol}
\begin{proof}
It is enough to check that if $\fE$ is H-nef and $c_1(E) \equiv
0$, then all universal quotient bundles $Q_{s,\fE}$ are nef.
Indeed, in this case the classes $\theta_{s,\fE} =
[c_1(Q_{s,\fE})]$ are nef, so that the classes $\lambda_{s,\fE}$
are nef by Theorem \ref{1.3BHR}. But
$\lambda_{s,\fE}=[c_1(\cO_{\PP Q_{s,\fE}}(1))]$, so that
$Q_{s,\fE}$ is nef.
\end{proof}
\begin{corol} Let $\fE=(E,\phi)$ be a  Higgs bundle on a smooth projective curve $X$,
such that all universal quotients $Q_{s,\fE}$ and $Q_{s,\fE^\ast}$ are nef.
Assume that $\fE$ is properly semistable, i.e., it is not stable. Then $\fE$ is
H-nflat.
\end{corol}

Proposition \ref{ss} raises the question of the existence
of stable H-nflat Higgs bundles. An example  bundle is provided by a Higgs bundle $\fE$  as in Example \ref{patho}
with $\deg(L_1)=1$, $\deg(L_2)=-1$ and the genus of the curve
$X$ at least 2. Then $\fE$ is stable and H-nflat.

Proposition \ref{ss} has a simple consequence, which generalizes
\cite[Thm. 1.18]{DPS94}. This provides an important
characterization of H-nflat Higgs bundles.

\begin{thm} A Higgs bundle $\fE$ on $X$ is H-nflat
if and only if it admits a filtration whose quotients are flat stable Higgs bundles.
\end{thm}
\begin{proof} If $\fE$ is H-nflat by Corollary \ref{corss}
it is semistable. Since all Chern classes of $\fE$ vanish, by     Theorem \ref{thm2} $\fE$
has a filtration whose quotients are stable and have vanishing Chern classes.
We may assume that $\fE$ is an extension
\begin{equation}\label{ext}  0 \to \fF \to \fE \to \fG \to 0 \end{equation} of stable
Higgs bundles with vanishing Chern classes, otherwise one simply iterates the following argument. Let us consider the bundle
$\fF=(F,\phi_F)$; the same will apply to $\fG$.  By results given in  \cite{S92}, the
bundle $F$ admits a Hermitian-Yang-Mills metric. Let $\Omega$ be the curvature of the
associated Chern connection. Since $c_1(F)=c_2(F)=0$, we have
$$ 0 = \int_X \mbox{tr}(\Omega\wedge \Omega) \cdot H^{n-2} =
\gamma_1 \Vert \Omega\Vert ^2 - \gamma_2\Vert\Lambda \Omega_i\Vert^2 = \gamma_1 \Vert
\Omega\Vert ^2$$ for some positive constants $\gamma_1$, $\gamma_2$, so that the Chern connection of $F$ is flat, i.e., $F$ is flat.

Conversely, let assume that $\fE$ has a filtration as in the statement. Then $\fE$ is
semistable with vanishing Chern classes, and by Corollary \ref{corral} it is
numerically flat.
\end{proof}

\end{document}